\documentclass{amsart}
\usepackage{amsmath,amssymb}
\usepackage[all]{xy}

\newtheorem{theorem}{Theorem}
\newcommand{\bt}{\begin{theorem}}
\newcommand{\et}{\end{theorem}}
\newtheorem{lemma}{Lemma}
\newcommand{\bl}{\begin{lemma}}
\newcommand{\el}{\end{lemma}}
\newtheorem{corollary}{Corollary}
\newcommand{\bc}{\begin{corollary}}
\newcommand{\ec}{\end{corollary}}
\newcommand{\beq}{\begin{equation}}
\newcommand{\eeq}{\end{equation}}
\newcommand{\benum}{\begin{enumerate}}
\newcommand{\eenum}{\end{enumerate}}
\newcommand{\N}{\ensuremath{ \mathbf N }}
\newcommand{\Z}{\ensuremath{\mathbf Z}}
\newcommand{\Q}{\ensuremath{\mathbf Q}}

\newcommand{\mcf}{\ensuremath{ \mathcal F}}

\newcommand{\mco}{\ensuremath{ \mathcal O}}

\newcommand{\mct}{\ensuremath{ \mathcal T}}
\DeclareMathOperator{\height}{\text{ht}}
\DeclareMathOperator{\id}{id}
\newcommand{\bmat}{\left(\begin{matrix}}
\newcommand{\emat}{\end{matrix}\right)}
\DeclareMathOperator{\qand}{\quad\text{and}\quad}
\DeclareMathOperator{\qqand}{\qquad\text{and}\qquad}

\title{A forest of linear fractional transformations}
\author{Melvyn B. Nathanson}
\address{Department of Mathematics\\
Lehman College (CUNY)\\
Bronx, NY 10468} 
\email{melvyn.nathanson@lehman.cuny.edu}

\subjclass[2010]{Primary 11A55, 11B75, 05A19, 05C05, 20M99.} 
\keywords{Calkin-Wilf tree, linear fractional transformation, 
continued fraction, $SL_2(\Z)$.}
\thanks{Supported in part by a grant from the PSC-CUNY Research Award Program.}

\date{\today}

\begin{document}

\begin{abstract}
The Calkin-Wilf tree is an infinite  binary tree whose vertices 
are the positive rational numbers.  
Each number occurs in the tree exactly once and in the form $a/b$, 
where are $a$ and $b$ are relatively prime positive integers.  
For every $2\times 2$ matrix  with nonnegative integral coordinates 
and nonzero determinant, it is possible to construct an analogous tree with this root.
If the root is the identity matrix, then the tree consists all matrices with determinant 1, 
and this tree possesses the basic properties 
of the Calkin-Wilf tree of positive rational numbers.  
The set of all matrices with nonzero determinant decomposes into a forest of rooted
infinite binary trees.  
\end{abstract}

\maketitle

\section{The Calkin-Wilf tree of rational numbers}
A \emph{directed graph}  consists of a nonempty set of \emph{vertices} 
and a set of directed \emph{edges}.  
Every edge is an ordered pair $(v,v')$ of vertices;  $v$ is called the \emph{tail} 
of the edge and $v'$ is called the \emph{head} of the edge.  
A \emph{path} in the graph from vertex $v$ to vertex $v'$ 
is a finite sequence of vertices $v=v_0, v_1,\ldots, v_k = v'$ 
such that, for all $i=1,\ldots, k$, either $(v_{i-1},v_i)$ is an edge 
or $(v_i, v_{i-1})$ is an edge.  The path is \emph{simple} if the vertices 
$v_0, v_1,\ldots, v_k$ are distinct.  
Distinct vertices $v$ and $v'$ are \emph{connected} if there is a simple path from $v$ to $v'$.
Every vertex is also connected to itself.  
This defines an equivalence relation on the vertices of the graph: 
$v$ and $v'$ are related if they are connected.
The equivalence classes of this relation are called the  
\emph{connected components} of the graph.  
A \emph{connected graph} is a graph in which every pair of edges is connected, 
that is, the set of all vertices is a component.
A directed graph is a \emph{tree} if there is a unique simple path between every pair of distinct vertices.

A \emph{rooted infinite  binary tree} is a tree with the following properties:
\benum
\item[(i)]
Every vertex is the tail of exactly two edges. 
Equivalently, every vertex has \emph{outdegree} 2.
\item[(ii)]
There is a vertex $v^*$ such that every vertex $v\neq v^*$ is the head 
of exactly one edge, but $v^*$ is not the head of any edge.  
Equivalently, every vertex $v\neq v^*$ has \emph{indegree} 1, and $v^*$ has indegree 0.
We call $v^*$ the \emph{root} of the tree.
\item[(iii)]
The graph is connected.
\eenum
A \emph{forest} is a directed graph whose connected components are 
rooted infinite  binary trees.

We call the rational number $a/b$ \emph{reduced} if $b\geq 1$ and 
the integers $a$ and $b$ are relatively prime.
The Calkin-Wilf tree~\cite{calk-wilf00} is an infinite  binary tree whose vertex set 
is the set of positive reduced rational numbers, 
and whose root is 1.  
In this tree, every positive reduced rational number 
$a/b$ is the tail of two edges.  
The heads of these edges are the positive rational numbers $a/(a+b)$ and $(a+b)/b$.  
Note that $a/(a+b) < 1 < (a+b)/b$.   We draw this as follows:
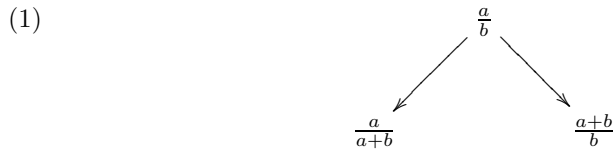
\beq   \label{Forest:parent-a/b}
\xymatrix{
& \frac{a}{b} \ar[dl]  \ar[dr] & \\
\frac{a}{a+b} &  & \frac{a+b}{b} \\
}
\eeq
with $a/(a+b)$ on the left and $(a+b)/b$ on the right.  
The fraction $a/b$ is the \emph{parent}.  
We call $a/(a+b)$ the \emph{left child} 
and $(a+b)/b$ the \emph{right child} of $a/b$.
Equivalently, if $z = a/b$, then the left child of $z$ is $z/(z+1)$ 
and the right child of $z$ is $z+1$.
These children give birth to children, and so on.  
Thus, every positive reduced  rational number has infinitely many descendants.  
If $\gcd(a,b)=1$, then $\gcd(a,a+b) = \gcd(a+b,b) = 1$,  
and so every descendant of a reduced rational number is also reduced. 
Equivalently, every positive reduced rational number is the root 
of an infinite  binary tree of positive reduced rational numbers.
The two children with the same parent are called \emph{siblings}.
The only positive rational number with no parent is 1, that is, 1 is an \emph{orphan}.  
Calkin and Wilf~\cite{calk-wilf00} introduced this enumeration 
of the positive rationals in 2000.  
It is related to the Stern-Brocot sequence~\cite{broc60,ster58}, 
discussed in~\cite{gram-knut-pata94}, 
and has stimulated much recent 
research (e.g. 
\cite{adam10, bate-bund-togn10, bate-mans11,dilc-stol07,glas11,mall11,mans-shat11,rezn90}).

The first four rows the Calkin-Wilf tree are as follows:
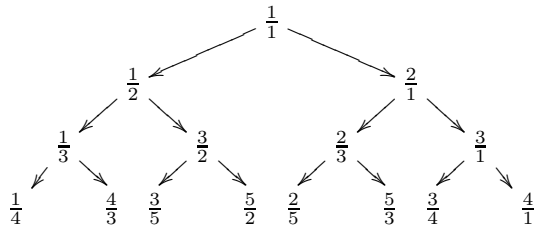
\[
\xymatrix@=0.2cm{
& & & & \frac{1}{1} \ar[dll]  \ar[drr] & & & & \\
& & \frac{1}{2} \ar[dl]  \ar[dr]  &  & & & \frac{2}{1} \ar[dl]  \ar[dr]  & & \\
& \frac{1}{3}   \ar[dl]  \ar[dr]  & & \frac{3}{2}  \ar[dl]  \ar[dr] &  & \frac{2}{3}  \ar[dl]  \ar[dr] &  & \frac{3}{1}  \ar[dl]  \ar[dr] & \\
\frac{1}{4}&  & \frac{4}{3} \quad \frac{3}{5} &  & \frac{5}{2} \quad  \frac{2}{5} &  & 
\frac{5}{3} \quad  \frac{3}{4} & & \frac{4}{1}  \\
}
\]
We enumerate the rows of the Calkin-Wilf tree as follows.
Row 0 contains only the number 1.  
Row 1 contains  the numbers 1/2 and 2.
For every nonnegative integer $n$, the $n$th row of the Calkin-Wilf tree 
contains $2^n$ positive reduced rational numbers. 
The $n$th row of the tree is also called the \emph{$n$th generation} of the tree.
We say that the rational number $a/b$ has \emph{depth} $n$, 
or belongs to generation $n$,  
if it is on the $n$th row of the tree.
We denote the ordered sequence of elements of the $n$th row, 
from left to right, by
$(c_{n,1}, c_{n,2}, \ldots, c_{n,2^n})$.  For example, 
$c_{2,3} = 2/3$ and $c_{3,6} = 5/3$.

Define the \emph{height} \index{height} of a nonzero rational 
number $x = a/b$ with $a$ and $b$ relatively prime integers 
 by $\height(x) = \max(|a|, |b|)$.  If $a/b$ is a positive reduced rational number, 
 then 
 \[
 \height\left(\frac{a}{a+b}\right) =  \height\left(\frac{a+b}{b}\right) = a+b \geq 
 \max(a,b) + 1 = \height\left(\frac{a}{b}\right) +1
 \]
and so the height of a child is strictly greater than the height of the parent.

\bt
Every positive reduced rational number appears exactly once 
as a vertex in the Calkin-Wilf tree.  
\et

\begin{proof}
The proof is by induction on the height.  
The only positive rational number of height 1 is 1.  
Because the height of a child is  strictly greater than the height of its parent, 
it follows that 1 is an orphan and appears only as the root of the Calkin-Wilf tree.  

Let $h \geq 2$, and suppose that every positive rational number of height less 
than $h$ appears exactly once as  a vertex in the tree.  
Let $a$ and $b$ be relatively prime positive integers, and let 
$a/b$ be a positive rational number of height $h = \max(a,b)$. 
Because $\max(a,b) \geq 2$, it follows that $a/b \neq 1$.
If $a/b < 1$, then $a < b$ and so $b-a \geq 1$.  
Thus,  $a/(b-a)$ is a positive reduced  rational number 
of height $\max(a,b-a) \leq b-1 < b = h$,  
and so $a/(b-a)$ is a vertex in the tree.
Moreover, $a/b$ is the left child of $a/(b-a)$, 
and so $a/b$ is a vertex in the tree.  
If $a'/b'$ is a fraction in the Calkin-Wilf tree 
such that $a/b$ is a child of $a'/b'$, then the inequality $a/b < 1$ 
implies that $a/b$ is the left child of $a'/b'$, and so 
\[
\frac{a}{b} = \frac{a'}{a'+b'}.
\]
Because $\gcd(a', a'+b')= 1$, it follows that $a'=a$ and $b'=b-a$.  
Thus, $a/b$ appears exactly once in the tree.

Similarly, if $a/b >1$, then $a-b \geq 1$ and $a/b$ is the right child 
of $(a-b)/b$, which is a  positive reduced rational number of height 
$\max(a-b,b) \leq a-1 < a = h$, 
and again $a/b$ appears exactly once as a vertex in the tree.  
This completes the proof.
\end{proof}

Here are four properties of the Calkin-Wilf tree:
\benum
\item
Denominator-numerator formula: 
For every positive integer $n$, we have $c_{n,1} = 1/(n+1)$ 
and  $c_{n,2^n} = n+1$.   
For $j=1,\ldots, 2^n-1$, if $c_{n,j} = p/q$, then $c_{n,j+1} = q/r$.  
Thus, as we move through the Calkin-Wilf tree from row to row, 
and from left to right across each row, the denominator of each fraction 
in the tree is the numerator of the next fraction in the tree.
This is in Calkin-Wilf~\cite{calk-wilf00}.

\item
Symmetry formula: 
For every nonnegative integer $n$ and for $j=1,\ldots, 2^n$, 
we have 
\[
 \frac{1}{ c_{n,j} } =  c_{n,2^n -j+1}.
\]
This straightforward observation is easily proved by induction on $n$.  

\item
Successor formula: 
For every positive integer $n$ and for $j=1,\ldots, 2^n -1$, we have
\[
c_{n,j+1} = \frac{1}{    2[c_{n,j}]+1-c_{n,j}  }
\]
where $[x]$ denotes the integer part of the real number $x$.  
This result is due to Moshe Newman~\cite{aign-zieg04,newm03}.

\item
Depth formula:
Let $a/b$ be a positive reduced rational number.
If 
\begin{align*}
\frac{a}{b} 
& = a_0 + \cfrac{1}{a_1 + \cfrac{1}{ a_2 + \cdots + \cfrac{1}{ a_{k-1}
+ \cfrac{1}{ a_k}}}} \\
& = [a_0, a_1,\ldots, a_{k-1}, a_k ]
\end{align*}
is the finite continued fraction of $a/b$, 
then the depth of $a/b$ is $a_0 + a_1 + \cdots + a_{k-1}  + a_k - 1$.
This is discussed in Gibbons, Lester, and Bird~\cite{gibb-lest-bird06}.
\eenum

For example, we have the continued fraction 
\[
\frac{11}{3} = 3 + \cfrac{1}{1+\cfrac{1}{2}} = [3,1,2]
\]
and so $11/3$ is on row $3+1+2 -1 = 5$ of the Calkin-Wilf tree.  
Because the integer part of $11/3$ is 3, 
the successor formula implies that the next element on row 5 is 
\[
\cfrac{1}{ 2\cdot 3+1- \cfrac{11}{3}} = \frac{3}{10}.
\]
Indeed, $c_{5,24} = 11/3$ and $c_{5,25} = 3/10$.  
Moreover, by the symmetry formula, $ 3/11 = 1/c_{5,24} = c_{5,9}$.

In this paper we describe a forest of rooted infinite  binary trees of rational functions 
of the form $(az+b)/(cz+d)$ in which all of these properties hold, 
and which specializes to the Calkin-Wilf tree 
when the root of the tree is the rational function $z=1$.

\section{Positive linear fractional transformations}
We begin with some definitions.  
Let $\N_0 = \{ 0,1,2,\ldots \}$ denote the set of nonnegative integers.  
Let $z$ be a variable.  
A \emph{positive linear fractional transformation} is a rational function 
\beq          \label{Forest:PLFT}     
f(z) = \frac{az+b}{cz+d}
\eeq
where $a,b,c,d$ are nonnegative integers 
such that $ad-bc \neq 0$.  
We call $ad-bc$ the \emph{determinant} 
\index{determinant!linear fractional transformation}
\index{linear fractional transformation!determinant} 
of $f(z)$, denoted $\det(f(z))$.  
Note that if $\det(f(z)) \neq 0$, then $(a,b) \neq (0,0)$ and $(c,d) \neq (0,0)$.  

Let 
\beq   \label{Forest:lft}  
f(z) = \frac{az+b}{cz+d} \quad \text{ and } \quad g(z) = \frac{ez+f}{gz+h}
\eeq
be positive linear fractional transformations, 
and consider the composite function
\[
f\circ g(z) 
= \frac{a \left(  \frac{ez+f}{gz+h}  \right) +b}{c  \left(  \frac{ez+f}{gz+h}  \right) +d} 
 = \frac{(ae+bg)z+(af+bh)}{(ce+dg)z+(cf+dh)}.
\]
The determinant of this composite function is 
\begin{align*}
\det(f\circ g(z)) 
& = (ae+bg)(cf+dh)  - (af+bh)(ce+dg) \\ 
& = (ad-bc)(eh-fg) \\
& = \det(f(z))\det(g(z)).
\end{align*}
Thus, the composite of positive linear fractional transformations 
is a positive linear fractional transformation.

A \emph{special positive linear  fractional transformation} 
is a positive linear fractional transformation 
$f(z)$ with $\det(f(z))=1$.
For example, 
\beq                    \label{Forest:l1r1}
\ell_1(z) = \frac{z}{z+1} \qqand 
r_1(z) = z+1
\eeq
are special positive linear fractional transformations.  
With the binary operation of composition, 
the set of positive linear fractional transformations is a monoid,
and the set of special positive linear fractional transformations 
is a submonoid of this monoid.

\bl        \label{Forest:lemma:SLT}
If $f(z)$ is a positive linear fractional  transformation, then $f(z)+1$ and 
$f(z)/(f(z)+1)$ are positive linear fractional  transformations, and
\[
\det(f(z)) = \det\left( f(z)+1  \right) = \det\left( \frac{f(z)}{f(z)+1}  \right).
\]
If $f(z)$ is a special positive linear fractional  transformation, then $f(z)+1$ and 
$f(z)/(f(z)+1)$ are special positive linear fractional  transformations.
\el

\begin{proof}
Let $f(z) = (az+b)/(cz+d)$ be a positive linear fractional transformation. 
The integers  $a,b,c,d$ are nonnegative with $\max(c,d) \geq 1$, and so
\[
f(z)+1 =  \frac{(a+c)z+ (b+d)}{cz+d} 
\]
is a positive linear fractional transformation with 
\[
\det(f(z)+1) = (a+c)d-(b+d)c = ad-bc = \det(f(z)).
\]
Similarly,  $\max(a,b) \geq 1$ and 
\[
\frac{f(z)}{f(z)+1} = \frac{az+b}{(a+c)z+ (b+d)} 
\]
is a positive linear fractional transformation with 
\[
\det\left( \frac{f(z)}{f(z)+1}   \right) = a(b+d) - b(a+c)
= ad-bc = \det(f(z)).
\]
This completes the proof.
\end{proof}

Lemma~\ref{Forest:lemma:SLT} also follows from the multiplicativity of the determinant 
and the observation that $f(z)+1 = r_1\circ f(z)$ and $f(z)/(f(z)+1) = \ell_1\circ f(z)$.

Let
\[
GL_2(\N_0) = \left\{
\bmat
a & b \\
c & d 
\emat : a,b,c,d \in \N_0 \text{ and } ad - bc \neq 0  
\right\}
\]
and
\[
SL_2(\N_0) = \left\{
\bmat
a & b \\
c & d 
\emat : a,b,c,d \in \N_0 \text{ and } ad - bc = 1 
\right\}.
\]
For example, 
\beq                    \label{Forest:L1R1}
L_1 = \bmat 1 & 0 \\1 & 1\emat \in SL_2(\N_0)
\qand
R_1 = \bmat 1 & 1 \\0 & 1\emat\in SL_2(\N_0).
\eeq
The sets $GL_2(\N_0)$ and $SL_2(\N_0)$ are monoids with the binary operation of matrix multiplication.  

The function that associates to every positive linear fractional 
transformation $f(z) =  (az+b)/(cz+d)$ the matrix
$\bmat a & b \\ c & d \emat \in GL_2(\N_0)$ is an isomorphism from 
the monoid of  positive linear fractional 
transformations onto $GL_2(\N_0)$.  
 The restriction of this function  to the monoid of special positive linear fractional 
transformation is an isomorphism onto $SL_2(\N_0)$.  

For example, the matrix associated to $f(z) = z$ is the identity matrix 
$I  = \bmat 1 & 0 \\ 0 & 1  \emat$
and the matrix associated to $f(z) = 1/z$ is the matrix 
$J  = \bmat 0 & 1 \\ 1 & 0  \emat$.  
If $f(z) =  (az+b)/(cz+d)$, then the matrices associated with 
$f(z)/(f(z)+1)$ and $f(z)+1$ are 
\[
L_1  \left( \begin{matrix}a & b \\c & d\end{matrix}\right)  
= \left( \begin{matrix}1 & 0 \\ 1 & 1 \end{matrix}\right)
\left( \begin{matrix}a & b \\c & d\end{matrix}\right) 
= \left( \begin{matrix}a & b \\ a+c & b+d \end{matrix}\right),
\]
and
\[
R_1 \left( \begin{matrix}a & b \\c & d\end{matrix}\right)  
= \left( \begin{matrix}1 & 1 \\0 & 1\end{matrix}\right) 
\left( \begin{matrix}a & b \\c & d\end{matrix}\right) 
= \left( \begin{matrix}a + c & b + d\\ c & d \end{matrix}\right)
\]
respectively.

\section{The tree of special positive linear fractional transformations}
\label{Forest:section:tree}
Let $z$ be a variable.  
We shall construct a rooted infinite  binary tree $\mct(z)$ with root $z$ 
whose vertex set is the monoid of special positive linear fractional transformations. 
Every vertex $w$ in this tree will be the parent of two children: 
the \emph{left child} $w /(w +1)$ 
and the \emph{right child} $w +1$.  
We draw this as follows:
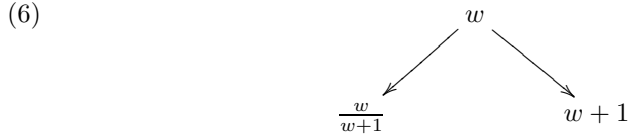
\beq   \label{Forest:parent}
\xymatrix{
& w  \ar[dl]  \ar[dr] & \\
\frac{w }{w +1} &  & w +1 \\
}
\eeq
with $w /(w +1)$ on the left and $w +1$ on the right.  
Note that if $w=a/b$ is a positive reduced rational number, 
then this is exactly the generation rule~\eqref{Forest:parent-a/b}.

The first four rows of $\mct(z)$ are as follows:
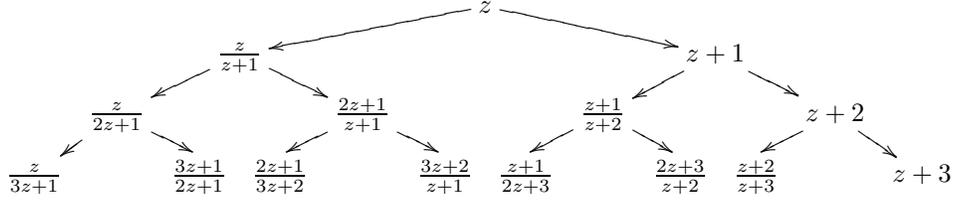
\[ 
\xymatrix@=0.16cm{
& & & & z  \ar[dll]  \ar[drr] & & & & \\
& & \frac{z }{z +1} \ar[dl]  \ar[dr]  &  & & & z +1 \ar[dl]  \ar[dr]  & & \\
&\frac{z }{2z +1}    \ar[dl]  \ar[dr]  & & \frac{2z +1}{z +1}  \ar[dl]  \ar[dr] &  &  \frac{z +1}{z +2}  \ar[dl]  \ar[dr] &  &z +2 \ar[dl]  \ar[dr] & \\
\frac{z }{3z +1} & &  \frac{3z +1}{2z +1} \quad \frac{2z +1}{3z +2} &  & \frac{3z +2}{z +1} \quad  \frac{z +1}{2z +3} &  & 
\frac{2z +3}{z +2} \quad  \frac{z +2}{z +3} & & z +3 
}
\]  
Because $z $ is a special linear fractional transformation, 
Lemma~\ref{Forest:lemma:SLT} implies 
that every vertex in this graph is a special linear fractional  transformation.
  
The root $z $ is the only element in row 0 of this tree.  
For every positive integer $n$, row $n$ of the tree consists 
of the $2^n$ elements of the $n$th generation descended 
from the root.  
We say that the rational function $f(z )$ has \emph{depth} $n$ 
if it is on row $n$ of the tree, or, equivalently, 
if it is a member of the $n$th generation of descendants 
of the root $z$.
We denote the ordered sequence of elements of the $n$th row by
$(w_{n,1}(z), w_{n,2}(z), \ldots, w_{n,2^n}(z))$.  For example, 
$w_{2,3}(z) = (z+1)/(z+2)$ and $w_{3,6}(z) =  (2z+3)/(z+2)$.
Note that $w_{2,3}(1) = 2/3 = c_{2,3}$ and $w_{3,6}(1) =  
5/3 = c_{3,6}$.
In general, $w_{n,j}(1) = c_{n,j}$ for all $n \in \N_0$ and $j \in \{1,2,\ldots, 2^n\}$.

A linear function $az+b$ will be  called  \emph{positive} 
if $a$ and $b$ are nonnegative integers and $(a,b) \neq (0,0)$.  
Equivalently, $az+b$ is positive if $a,b \in \N_0$ and $a+b > 0$.   
We partially order the set of positive linear functions  as follows:
\[
cz+d \preceq az+b \qquad \text{if $c \leq a$ and $d \leq b$.}
\]
We write $cz+d \prec az+b$ if $cz+d \preceq az+b$  
and $cz+d  \neq az+b$.
Distinct positive linear functions $az+b$ and $cz+d$ 
are \emph{comparable} if $cz+d \prec az+b$ or $az+b \prec cz+d$.
For example, $2z+1$ and $3z+2$ are comparable, 
but $2z+1$ and $z+2$ are not comparable.  
A positive linear fractional transformation $(az+b)/(cz+d)$ 
such that $az+b \prec cz+d$ is the left child 
of the positive linear transformation $(az+b)/( (c-a)z+(d-b))$.  
A positive linear fractional transformation $(az+b)/(cz+d)$ 
such that $cz+d \prec az+b$ is the right child 
of the positive linear transformation 
$((a-c)z+(b-d))/(cz + d)$.  
If $az+b$ and $cz+d$ are not comparable, 
then the positive linear fractional transformation 
$(az+b)/(cz+d)$ has no parent and  is called an \emph{orphan}.  
For example, $z$, $(z+2)/(2z+1)$ and $1/3z$ are orphans 
with determinants 1, -3, and -3, respectively.  

We define the \emph{height} of the positive linear fractional transformation 
$f(z) = (az+b)/(cz+d)$ by $\height(f(z)) = \max(a+b,c+d)$.
If $g(z)$ is the left or right child of $f(z)$, then $\height(g(z)) = a+b+c+d > \height(f(z))$.
The height of a positive linear fractional transformation 
is a positive integer, and the height of a child is always strictly greater 
than the height of its parent.

\bt     \label{Forest:theorem:SLFTtree}
The directed graph $\mct(z)$ is a rooted infinite binary tree 
with root $z$.  The set of vertices of the tree $\mct(z)$ 
is the set of all special positive linear fractional transformations, 
and each special positive linear fractional transformation 
occurs exactly once as a vertex in this tree.
\et

\begin{proof}
We have already observed that every vertex of $\mct(z)$ is a special positive linear fractional transformation.  Conversely, every special positive linear transformation 
$f(z) = (az+b)/(cz+d)$  is either an orphan or has a parent.  
That parent is either an orphan or has a parent.  
Because the height of a parent is always strictly less than the height of its children, 
it follows that every vertex in the tree $\mct(z)$ has only finitely many ancestors, 
and so every vertex is the descendent of an orphan.  

If $f(z)$ is an orphan, then the positive linear functions $az+b$ and $cz+d$ are not comparable.  This can happen in only two ways.  
In the first case, we have $a<c$ and $b > d$, and so
\[
1 = \det(f(z)) = ad-bc \leq (c-1)(b-1) - bc = 1 - b - c \leq 1
\] 
and so $b=c=0$ and $a=d=1$, hence $f(z)=z$.

In the second case, we have $a > c$ and $b < d$, and so 
\[
1 = ad-bc \geq (c+1)(b+1) - bc = 1 + b + c \geq 1
\]
and so $b=c=0$ and $a=d=1$, hence $f(z)=z$.
Thus, the only orphan special positive linear transformation is $z$.  
Because every positive linear fractional transformation is descended 
from an orphan, every special positive linear fractional transformation 
is a descendent of $z$, and must be a vertex in the tree $\mct(z)$.  
Moreover, every vertex has a unique parent, and so every 
special positive linear fractional transformation occurs exactly once 
as a vertex in the tree $\mct(z)$.   
This completes the proof.
\end{proof}

\bc
The monoid $SL_2(\N_0)$ is a free monoid of rank 2, 
and the set of matrices $\{L_1,R_1\}$ defined by~\eqref{Forest:L1R1}
freely generates $SL_2(\N_0)$.
\ec

\begin{proof}
This follows immediately from the isomorphism between the monoid of  
special positive fraction linear transformations and the fact 
(Theorem~\ref{Forest:theorem:SLFTtree}) 
that the graph $\mct(z)$ is a tree.  
\end{proof}

\section{Properties of the tree $\mct(z)$}

We shall prove that, with appropriate definitions 
of ``integer part,'' ``reciprocal,'' and ``continued fraction,'' 
properties~(1)-(4) of the Calkin-Wilf positive rational number tree also hold 
for the tree $\mct(z)$ of special positive linear fractional transformations.
Recall that, for $j = 1,\ldots, 2^n$, the linear fractional 
transformation $w_{n,j}(z)$ is the $j$th vertex on the $n$th row of $\mct(z)$.  

\bt[Denominator-numerator formula]
For all $n \geq 1$ and $j=1,\ldots, 2^n -1$, 
the denominator of $w_{n,j}(z)$ is the numerator of $w_{n,j+1}(z)$.
\et

\begin{proof}
The proof is by induction on $n$.  
The theorem is true for $n=1$ because $z+1$ is both the denominator of $w_{1,1}(z)$ 
and the numerator of $w_{1,2}(z)$.  

Let $n\geq 2$, and assume that the theorem holds for $n-1$.  
If $j$ is odd, then $w_{n,j}(z)$ and $w_{n,j+1}(z)$  are siblings.  
If their parent is the linear fractional transformation $(az+b)/(cz+d)$, 
then $(a+c)z+(b+d)$ is  the denominator of 
$w_{n,j}(z)$ and the numerator of $w_{n,j+1}(z)$.  

If $j=2i$ is even, then $i$ is a positive integer with $i < 2^{n-1}$ 
such that 
$w_{n,j}(z)$ is the right child of $w_{n-1,i}(z)$ and $w_{n,j+1}(z)$ is the 
left child of $w_{n-1,i+1}(z)$.  
If $w_{n-1,i}(z) = (az+b)/(cz+d)$, then the induction hypothesis 
implies that $w_{n-1,i+1}(z) = (cz+d)/(ez+f)$.
The right child of $w_{n-1,i}(z)$  is 
$w_{n,j}(z) = ((a+c)z+(b+d))/(cz+d)$; 
the left child of $w_{n-1,i+1}(z)$  is 
$w_{n,j+1}(z) = (cz+d)/((c+e)z+(d+f))$.  
We see that $cz+d$ is both the denominator of 
$w_{n,j}$ and the numerator of $w_{n,j+1}(z)$.  
This completes the proof.  
\end{proof}

\bt[Symmetry formula]
Define the function $\Phi$ on the set of nonzero rational functions in $z$ by 
\[
\Phi(f(z)) = \frac{1}{f\left(\frac{1}{z}\right)}.
\]
Then $\Phi$ is an involution, that is, $\Phi^2 = \id$, and, 
for every positive integer $n$ and $j=1,\ldots, 2^n$,
\beq   \label{Forest:symPhi}
\Phi(w_{n,j})(z) = w_{n,2^n -j+1}(z).
\eeq
\et

\begin{proof}
For every rational function $f(z)$ we have 
\[
\Phi^2( f(z)) = \Phi\left( \Phi(f(z) )\right) = \frac{1}{\Phi(f(1/z))} 
= \frac{1}{\frac{1}{f\left(\frac{1}{1/z}\right)}} = f(z)
\]
and so $\Phi^2 = \id$.

We shall prove~\eqref{Forest:symPhi} by induction on $n$.
We have  
\[
\Phi\left( w_{1,1}(z) \right) = \Phi\left( \frac{z}{z+1} \right) = 
\frac{1}{ \frac{ \frac{1}{z} }{\frac{1}{z} +1} } = z+1 = w_{1,2}.
\]
Because $\Phi$ is an involution, we have $\Phi(w_{1,2}(z)) = w_{1,1}(z)$. Thus,~\eqref{Forest:symPhi} holds for $n=1$.

In general, if $f(z) = (az+b)/(cz+d)$, the $\Phi(f(z)) = (dz+c)/(bz+a)$.
Let $n \geq 2$ and suppose that~\eqref{Forest:symPhi} holds for $n-1$ and $j=1,\ldots, 2^{n-1}$.  
If 
\[
w_{n-1,j}(z) = \frac{az+b}{cz+d}
\]
then 
\[
w_{n-1,2^{n-1}-j+1}(z)  = \Phi\left(  w_{n-1,j}(z) \right) 
 = \Phi\left( \frac{az+b}{cz+d}\right) 
= \frac{dz+c}{bz+a}.
\]
The children of $w_{n-1,j}(z)$ are 
\[
w_{n,2j-1}(z) = \frac{az+b}{(a+c)z+(b+d)}
\]
and 
\[
w_{n,2j}(z)  = \frac{(a+c)z+(b+d)}{cz+d}.
\]
The children of $w_{n-1,2^{n-1}-j+1}(z)$ are 
\[
w_{n,2^{n}-2j+1}(z) =  \frac{dz+c}{(b+d)z+(a+c)}
\]
and 
\[
w_{n,2^{n}-2j+2}(z)  =  \frac{(b+d)z+(a+c)}{bz+a}.
\]
We see immediately that 
\[
\Phi\left( w_{n,2j-1}(z) \right) = w_{n,2^{n}-2j+2}(z) 
\]
and
\[
\Phi\left( w_{n,2j}(z)  \right) = w_{n,2^{n}-2j+1}(z).
\]
This completes the proof.  
\end{proof}

The \emph{greatest common divisor} of the positive linear functions
$az+b$ and $cz+d$ is their greatest common divisor in the ring of polynomials 
with integer coefficients.  
For example, the  greatest common divisor of $10z+6$ and $15z+9$ is $5z+3$ 
and the  greatest common divisor of $9z+6$ and $15z+9$ is 3.
The positive linear functions
$az+b$ and $cz+d$ are \emph{relatively prime} if their greatest common divisor is 1.

\bl[Division algorithm]    \label{Forest:lemma:DivisionAlgorithm}  
Let $az+b$ and $cz+d$ be relatively prime positive linear functions.  
If $cz+d \prec  az+b$, then there is a unique positive integer $q$ 
and  a unique positive linear function $rz+s$   
such that the polynomials $rz+s$ and $cz+d$ are relatively prime, 
\beq    \label{Forest:DivisionAlgorithm}
az+b = q(cz+d) + (rz+s)
\eeq
and  either $r < c$ or $s < d$.
Moreover, $rd-sc = ad-bc$.
\el

Equivalently,
\[
\frac{az+b}{cz+d} = q + \frac{rz+s}{cz+d}
\]
with $r < c$ or $s < d$, and 
\[
\det\left( \frac{az+b}{cz+d} \right) = \det\left(  \frac{rz+s}{cz+d} \right).
\]

\begin{proof}
Because $cz+d \prec az+b$ and $az+b$ is not a multiple of $cz+d$, 
there is a largest positive integer $q$ 
such that $ q(cz+d) \prec az+b$.  
Indeed, if $cd > 0$, then $q = \min([a/c],[b/d]) \geq 1$.
If $c=0$, then $q = [b/d] \geq 1$, and if $d=0$, then $q = [a/c] \geq 1$.  
In all three cases, we define $r = a-qc$ and $s = b - qd$.
Then $r$ and $s$ are nonnegative integers such that 
\[
rz+s = (a-qc)z + (b - qd) = (az+b) - q(cz+d).
\]
If $q=[a/c]$, then  $a/c < q+1$ and $r =a-qc < c$.  
If $q = [b/d]$, then $b/d < q+1$ and $s = b-qd <d$.
Thus, either $r < c$ or $s < d$.
Moreover,
\[
rd-sc = (a-qc)d - (b-qd)c = ad-bc  
\]
and so the  $(az+b)/(cz+d)$ 
and  $(rz+s)/(cz+d)$ have the same determinant.  

Let $q$ and $q'$ be a positive integers and let $rz+s$ and $r'z + s'$ 
be positive linear functions satisfying 
\[
az+b = q(cz+d) + rz+s = q'(cz+d) + r'z+s'
\]
with the property that $r < c$ or $s < d$, and also that  $r' < c$ or $s' < d$.
Then 
\[
qc+r = q'c+r'  \quad \text{ and } \quad qd+s = q'd+ s'.
\]
If $q' > q$, then 
\[
qc+r = q'c+r' \geq (q+1)c+r' = qc+c+r'
\]
and so $r \geq c+r' \geq c$.   Similarly,
\[
qd+s = q'd+s' \geq (q+1)d+s' = qd + d + s'
\]
and $s \geq d+s' \geq d$.  
Thus, if $q' > q$, then $r \geq c$ and $s \geq d$, which is absurd.
Similarly, if $q' < q$, then $r' \geq c$ and $s' \geq d$, 
which is also absurd.  
It follows that $q=q'$ and $rz+s = r'z+s'$.   
This completes the proof.  
\end{proof}

In the division algorithm~\eqref{Forest:DivisionAlgorithm}, 
we call $q$ the \emph{integer part} of the linear fractional transformation 
$(az+b)/(cz+d)$ and write
\[
q = \left[ \frac{az+b}{cz+d}\right].
\]
We call $(rz+s)/(cz+d)$ the \emph{fractional part} 
of the linear fractional transformation $(az+b)/(cz+d)$ and write
\[
\frac{rz+s}{cz+d} = \left\{ \frac{az+b}{cz+d}\right\}.
\]
For example, if
\[
f(z) = \frac{21z+16}{8z+5} 
\]
then the division algorithm gives 
\[
21z+16 = 2(8z+5) + 5z+6.
\]
The integer part of $f(z)$ is 2 and the fractional part of $f(z)$ 
is $(5z+6)/(8z+5)$.  Note that $8z+5 \prec 21z+16$, 
but $8z+5$ and $5z+6$ are not comparable, that is, $\{f(z) \}$ is an orphan.

If $az+b = q(cz+d)$, then we say that  the integer part of $(az+b)/(cz+d)$ is $q$ 
and the fractional part is $0$.
If $az+b \prec cz+d$, then we say that  the integer part of $(az+b)/(cz+d)$ is $0$ 
and the fractional part is $(az+b)/(cz+d)$.  
Note that the integer and fractional parts of $(az+b)/(cz+d)$ are 
undefined if $az+b$ and $cz+d$ are unequal and not comparable.

\bl   \label{Forest:lemma:iterate}
Let $w$ be a vertex in the infinite binary tree generated by $z$. 
The descendant of $w$ after $k$ generations to the right  is $w+k$, 
and the descendant of $w$ after $k$ generations  to the left is 
$w/(kw+1)$.
\el

\begin{proof}
For $k=1$ this is simply the definition of the right and left descendants.  

Let $k \geq 2$.  
If the right descendant of $w$ after $k-1$ generations is $w+k-1$, 
then the right descendant of $w$ after $k$ generations is $w+k$.  
If the left descendant of $w$ after $k-1$ generations is 
$w/((k-1)w+1)$, 
then the left descendant of $w$ after $k$ generations is 
\[
\frac{  \frac{w}{(k-1)w+1}}{ \frac{w}{(k-1)w+1}+1} 
= \frac{w}{kw+1}.
\] 
This completes the proof.  
\end{proof}

\bt[Successor formula]
Let $n$ be a positive integer.  
In the infinite binary tree generated by $z$, 
if $w_{n,j}(z)$ and $w_{n,j+1}(z)$ are successive terms 
on the $n$th row of the linear fractional transformation tree, then 
\[
w_{n,j+1}(z) = \frac{1}{2[w_{n,j}(z)] + 1 - w_{n,j}(z) }
\]
where $[w_{n,j}(z)] $ is the integer part of $w_{n,j}(z)$.  
\et

\begin{proof}
Let $i \in \left\{1,2,\ldots, 2^{n-1} \right\}$ and $j = 2i-1$. 
The linear fractional transformations $w_{n,2i-1}(z)$ and $w_{n,2i}(z)$ 
are successive elements on the $n$th row, 
and are the left and right children of $w_{n-1,i}(z)$.  
If $w_{n-1,i}(z) = (az+b)/(cz+d)$, then 
\[
w_{n,2i-1}(z) = \frac{az+b}{(a+c)z+(b+d)}.
\]
Because $az+b \prec (a+c)z+(b+d)$, we have 
$[w_{n,2i-1}(z)] = 0$ and $\{w_{n,2i-1}(z)\} = w_{n,2i-1}(z).$
Then
\begin{align*}
\frac{1}{2[w_{n,2i-1} (z)] + 1 - w_{n,2i-1}(z) } 
& =  \frac{1}{1 - w_{n,2i-1}}  \\
& = \frac{1}{1 -\frac{az+b}{(a+c)z+(b+d)}  } \\
& = \frac{(a+c)z+(b+d)}{cz+d} \\
& = w_{n,2i}(z).
\end{align*}

Let $i \in \left\{1,2,\ldots, 2^{n-1} -1 \right\}$ and $j = 2i$. 
If $w_{n,2i}(z)$ and $w_{n,2i+1}(z)$ are successive elements on the 
$n$th row, then the former is the right child and the latter is the left child 
of successive elements in the $(n-1)$st row.   
If these linear fractional transformations on the $(n-1)$st row are not siblings, 
then they are the right and left children, respectively, of successive elements in row $n-2$.  
Every element in the tree is a descendant of the root $z$.  
Retracing the family tree, we must eventually reach an element 
from which both $w_{n,2i}(z)$ and $w_{n,2i+1}(z)$ are descended.  
Thus, there is a smallest positive integer $k$ such that this common 
ancestor is on row $n-k-1$.  Let $w^* = w_{n-k-1,t}(z)$ be this ancestor.  
Its children are $w_{n-k,2t-1}(z)$ and $w_{n-k,2t}(z)$.
Then $w_{n,2i}(z)$ is the $k$-fold right child of $w_{n-k,2t-1}(z)$, 
and $w_{n,2i+1}(z)$ is the $k$-fold left child of $w_{n-k,2t}(z)$.
Thus, 
\[
w_{n-k,2t-1}(z) = \frac{w^*}{w^*+1} 
\]
and, by Lemma~\ref{Forest:lemma:iterate} with $w = w_{n-k,2t-1}$,
\begin{align*}
w_{n,2i}(z) & = w_{n-k,2t-1} + k 
= \frac{w^*}{w^*+1} + k.
\end{align*}
Because $w^* \prec w^*+1$,  the division algorithm 
(Lemma~\ref{Forest:lemma:DivisionAlgorithm})
implies that 
\[
[w_{n,2i}(z) ] = k.
\]
Similarly,
\[
w_{n-k,2t}(z) = w^* + 1
\]
and so, by Lemma~\ref{Forest:lemma:iterate}  with $w = w_{n-k,2t}$,
\begin{align*}
w_{n,2i+1}(z) 
& = \frac{w_{n-k,2t}(z) }{kw_{n-k,2t}(z) +1} \\
& = \frac{w^* + 1}{k(w^* + 1) +1} \\
& = \frac{1}{k+ 1 - \frac{w^*}{w^* + 1}} \\
& = \frac{1}{2k + 1 - w_{n,2i}(z) } \\
& = \frac{1}{ 2[w_{n,2i}(z) ]  + 1 - w_{n,2i}(z) }.
\end{align*}
This completes the proof.
\end{proof}

\section{Continued fractions and the depth formula}
To prove the analogue of the depth formula, we need to introduce 
finite continued fractions of linear fractional transformations.  

Let $az+b$ and $cz+d$ be comparable relatively prime 
positive linear functions, that is,  
either $cz+d \prec az+b$ or $az+b \prec cz+d$.  
Note that if $cz+d \prec az+b$, then $0 < c+d < a+b$.
We define
\[
r_0z+s_0 = az+b
\]
and
\[
r_1z+s_1 = cz+d.
\]
If $r_0z+s_0 \prec   r_1z+s_1$, then we have
\[
r_0z+s_0 = q_0 (r_1z+s_1) + (r_2z+s_2)
\]
where $q_0 = 0$, $r_2z+s_2 = r_0z+s_0 $, 
and $r_2z+s_2 \prec r_1z+s_1$.
If  $r_1z+s_1\prec r_0z+s_0$, 
then, by the division algorithm (Lemma~\ref{Forest:lemma:DivisionAlgorithm}),  there exist a unique positive integer $q_0$ 
and a unique positive linear function $r_2z+s_2$ such that 
\[
r_0z+s_0 = q_0 (r_1z+s_1) + (r_2z+s_2)
\]
and either $r_2z+s_2 \prec r_1z+s_1$ 
or the linear functions $r_1z+s_1$ and $r_2z+s_2 $ are not comparable.  

If $r_2z+s_2 \prec r_1z+s_1$, then, first,  $0 < r_2+s_2 < r_1+s_1$, 
and, second,  there exist a unique positive integer $q_1$ 
and a unique positive linear function 
$r_3z+s_3$ such that 
\[
r_1z+s_1 = q_1 (r_2z+s_2) + (r_3 z+s_3)
\]
and either $r_3z+s_3 \prec r_2z+s_2$ (and so $0< r_3+s_3 < r_2+s_2$) 
or the linear functions $r_2 z+s_2 $ and $r_3z+s_3 $ are not comparable.  

Continuing inductively, we obtain a finite sequence of positive 
linear functions $r_iz+s_i$ for $i=0,1,\ldots, k+1$ such that
$r_i z+s_i \prec r_{i-1}z+s_{i-1}$ for $i=2,3,\ldots, k$ and 
\begin{align*}
r_0z+s_0 & = q_0 (r_1z+s_1) + (r_2z+s_2) \\ 
r_1z+s_1 & = q_1 (r_2z+s_2) + (r_3 z+s_3) \\
r_2 z+s_2 & = q_2 (r_3 z+s_3 ) + (r_4  z+s_4 ) \\
& \vdots \\
r_{k-2} z+s_{k-2} & = q_{k-2} (r_{k-1} z+s_{k-1} ) + (r_{k}  z+s_{k} ) \\
r_{k-1} z+s_{k-1} & = q_{k-1} (r_{k} z+s_{k} ) + (r_{k+1}  z+s_{k+1} ). 
\end{align*}
The linear function $r_{k+1}  z+s_{k+1}$ is nonzero 
because the positive linear functions $az+b$ and $cz+d$ are relatively prime.  
Because every strictly decreasing sequence of positive integers is finite and 
\[
0 < r_k+s_k < r_{k-1}+s_{k-1} < \cdots < r_1+s_1 
\]
the process of iteration of the division algorithm  must terminate, 
and, after, say, $k$ divisions, we obtain positive linear forms 
$r_{k} z+s_{k} $ and $r_{k+1}  z+s_{k+1}$ that are not comparable.  
We call this procedure the \emph{Euclidean algorithm}.

We rewrite the equations in the Euclidean algorithm 
to obtain a \emph{finite continued fraction} for the linear fractional 
transformation
\begin{align*}
\frac{az+b}{cz+d} 
& = \frac{r_0z+s_0}{r_1 z+s_1} \\
& = q_0 + \frac{r_2z+s_2}{r_1z+s_1} 
 = q_0 + \cfrac{1}{\cfrac{r_1z+s_1}{r_2z+s_2}} \\
& = q_0 + \cfrac{1}{q_1 + \cfrac{r_2z+s_2}{r_3 z+s_3}} 
 = q_0 + \cfrac{1}{q_1 + \cfrac{1}{\cfrac{r_3 z+s_3}{r_2z+s_2}} }\\
& \vdots \\
& = q_0 + \cfrac{1}{q_1 + \cfrac{1}{  q_2 + \cdots + \cfrac{1}{ q_{k-1} 
+ \cfrac{1}{\cfrac{ r_k  z+s_k}{ r_{k+1}z+s_{k+1}  } } }  }}.
\end{align*}
The linear fractional transformation $w^* = ( r_k  z+s_k)/(r_{k+1}z+s_{k+1})$ 
is an orphan; we call it the \emph{root} of the  linear fractional transformation $(az+b)/(cz+d)$.  
We call $[q_0,q_1,\ldots, q_{k-1},w^*]$ the \emph{continued fraction} 
of $(az+b)/(cz+d)$ with root $w^*$.  
Note that  $q_0=0$ if $(az+b)/(cz+d)$ is a left child 
and $q_0\geq 1$ if $(az+b)/(cz+d)$ is a right child.

For example, applying the Euclidean algorithm to the relatively prime 
positive linear functions $21z+46$ and $5z+11$, we obtain 
\begin{align*}
21z+46 & = 4(5z+11) + (z +2)\\
5z+11 & = 5(z +2) +1
\end{align*}
and so
\[
f(z) = \frac{21z+46}{5z+11} = 4 + \frac{z+2}{5z+11} 
= 4 + \cfrac{1}{5+\cfrac{1}{2+z}}.
\]
Thus, the continued fraction of the positive linear fractional transformation 
$f(z)$ is $[4,5,2+z]$.

Consider the positive linear fractional transformation 
$g(z) = (17z+10)/(5z+3)$.
Applying the Euclidean algorithm to the relatively prime 
positive linear functions $17z+10$ and $5z+3$, we obtain 
\begin{align*}
17z+10 & = 3(5z+3) + (2z+1) \\
5z+3 & = 2(2z + 1) + (z+1) \\
2z + 1 & = 1( z+1) + z \\
z+1& = 1( z )+1
\end{align*}
and so
\begin{align*}
g(z) & = \frac{17z+10}{5z+3} = 3 + \frac{2z+1}{5z+3} 
= 3 + \cfrac{1}{2+\cfrac{z+1}{2z+1}} \\
& = 3 + \cfrac{1}{2+\cfrac{1}{1+\cfrac{z}{z+1}}}
= 3 + \cfrac{1}{2+\cfrac{1}{1+\cfrac{1}{1+\cfrac{1}{z}}}}.
\end{align*}
The continued fraction of 
$g(z)$ is $[3,2,1,1,z]$.

Revisit the first four rows ($n = 0,1,2,3$) of the tree with root $z$ from 
the beginning of Section~\ref{Forest:section:tree}.
We denote the elements of the $n$th row, from left to right, by 
$w_{n,1}(z), w_{n,2}(z), \ldots, w_{n,2^n}(z)$.
The continued fractions of the special positive 
linear fractional transformations in these rows are as follows:\\
\[
\begin{array}{ |    lll | l  lll     |}
\hline
w_{0,1}(z) & = z & = [z] & \quad 
& w_{3,1}(z)  & =   \cfrac{1}{3+\cfrac{1}{z}} 
& =  [0,3,z]\\
w_{1,1}(z)  &  = \cfrac{1}{1+\cfrac{1}{z}} & = [0,1,z]  & \quad 
& w_{3,2}(z)  & = 1 +  \cfrac{1}{2+\cfrac{1}{z}}  & =  [1,2,z]  \\
w_{1,2}(z)  & = 1+ z & = [1+z]  & \quad  
& w_{3,3}(z)  & =  \cfrac{1}{1+  \cfrac{1}{1+ \cfrac{1}{1+ \cfrac{1}{z}}}}
& = [0,1,1,1,z] \\
w_{2,1}(z)  &  = \cfrac{1}{2+\cfrac{1}{z}} 
& = [0,2,z] & \quad 
& w_{3,4}(z)  & =  2 + \cfrac{1}{1 + \cfrac{1}{z}} 
& = [2,1,z]   \\
w_{2,2}(z)  & = 1 + \cfrac{1}{1+\cfrac{1}{z}} 
& = [1,1,z] & \quad  
& w_{3,5}(z)  & = \cfrac{1}{2+\cfrac{1}{1+z}}
& = [0,2,1+z]   \\
w_{2,3}(z)  & =  \cfrac{1}{1+\cfrac{1}{1+z}}
& = [0,1,1+z] & \quad  
& w_{3,6}(z)  & =  1 + \cfrac{1}{1 + \cfrac{1}{1+z}}
& = [1,1,1+z]  \\
w_{2,4}(z)  & = 2 + z & = [2 + z]  & \quad 
& w_{3,7}(z)  & =   \cfrac{1}{1 + \cfrac{1}{2+z}} 
& =  [0,1,2+z]  \\
& & & \quad &  w_{3,8}(z)  & = 3+z & =  [3 + z]\\
&&&&&&\\
 \hline
\end{array}
\]
\\
We can write each vertex of the tree $\mct(z)$ as a finite continued fraction.  
The  first three rows of left descendants of the vertex are:
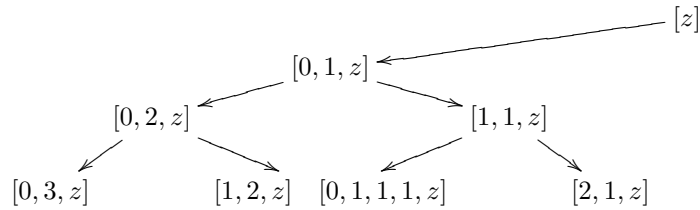
\[
\xymatrix@=0.1cm{
& & & & &[z]  \ar[dlll] \\
& & [0,1,z]\ar[dl]  \ar[dr]  &  & & \\
& [0,2,z]    \ar[ddl]  \ar[ddr]  & &  [1,1,z]  \ar[ddl]  \ar[ddr] & & \\
&&&&& \\
[0,3,z]  & & [1,2,z] \quad [0,1,1,1,z] &  &  [2,1,z] &
}
\]
The  first three rows of right descendants of the vertex are:
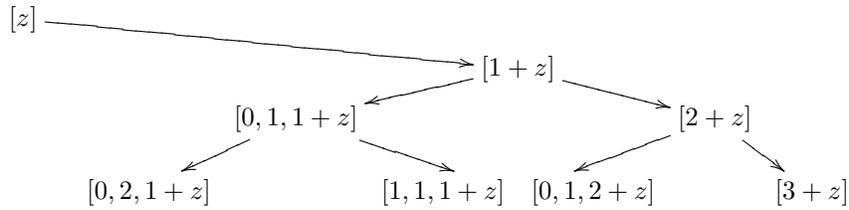
\[
\xymatrix@=0.1cm{
 [z]    \ar[drrr] & & & & &  \\
&  & &[1+z]  \ar[dl]  \ar[dr]  & & \\
& &  [0,1,1+z]  \ar[ddl]  \ar[ddr] &  & [2+z]\ar[ddl]  \ar[ddr] & \\
& & & & & \\
& \quad [0,2,1+z] &  & [1,1,1+z] \quad  [0,1,2+z] & & [3+z]
}
\]

\bt[Depth formula]
Every vertex $f(z)$ in  the infinite binary tree generated by $z$ 
has a unique continued fraction in exactly one of the following two forms:
\beq       \label{Forest:StandardFormA}
f(z) = [q_0, q_1,\ldots, q_k +  z]
\eeq
with $k$ even, or 
\beq       \label{Forest:StandardFormB}
f(z) = [q_0, q_1,\ldots, q_k,  z]
\eeq
with $k$ odd.
Moreover, $f(z)$ is in row $q_0 + q_1 + \cdots + q_k$ of the tree.  
\et

\begin{proof}
The unique element of row 0 is the root $z$, whose continued fraction 
$z =  [z]$ is of the form~\eqref{Forest:StandardFormA} 
with $k = 0$ and $q_0=0$.  
Similarly, the rational functions on row 1 are $z/(z+1) = [0,1,z]$ with $k=1$, and 
$z+1 = [1+z]$ with $k=0$.  

Let $n \geq 1$, and assume that the Theorem is true for the rational functions 
on the $n$th row.  Let $v  = f(z)$ be on row $n+1$.  If $v$ is a right child, then 
there exists $v'$ on row $n$ such that $v = v'+1$.
If $v' = [q'_0, q'_1,\ldots, q'_k +  z] $ is of form~\eqref{Forest:StandardFormA},
then $k$ is even, $\sum_{i=0}^k q'_i = n$, and 
\[
v = v'+1 = [q'_0, q'_1,\ldots, q'_k +  z] + 1 = [q_0, q_1,\ldots, q_k +  z]
\]
with $q_0 = q'_0+1$, and $q_i = q'_i$ for $i=1,\ldots, k$.  
Similarly, if $v =  [q_0, q_1,\ldots, q_k,  z]$ 
is of form~\eqref{Forest:StandardFormB}, then $k$ is odd, 
$\sum_{i=0}^k q'_i = n$, and 
\[
v = v' +1 = [q'_0, q'_1,\ldots, q'_k, z] + 1 = [q_0 , q_1,\ldots, q_k,  z]
\]
with $q_0 = q'_0+1$, and $q_i = q'_i$ for $i=1,\ldots, k$.  
In both cases, 
$\sum_{i=0}^k q_i = 1 + \sum_{i=0}^k q'_i = n+1$.

If $v$ is a left child on the $(n+1)$st row, 
then there exists $v'$ on row $n$ such that $v = v'/(v' +1)$.
Let $v' = [q'_0, q'_1,\ldots, q'_k +  z] $ be of form~\eqref{Forest:StandardFormA}, 
with $k$ is even and $\sum_{i=0}^k q'_i = n$.  
If $q'_0 \geq 1$, then 
\begin{align*}
v & = \frac{v'}{v'+1} = \cfrac{1}{1+\cfrac{1}{v'}} 
 = \cfrac{1}{1+\cfrac{1}{[q'_0, q'_1,\ldots, q'_k +  z] }} \\
& = [0,1,q'_0, q'_1,\ldots, q'_k +  z] \\
& =  [q_0, q_1,q_2, \ldots, q_{k+1},q_{k+2} +  z] 
\end{align*}
with $q_0 = 0$, $q_1 = 1$, and $q_i = q'_{i-2}$ for $i=2,3,\ldots, k+2$.  
Moreover, $\sum_{i=0}^{k+2} q_i = 1 + \sum_{i=0}^k q'_i = n+1$.

If $q'_0 =0$, then 
\begin{align*}
v & = \frac{v'}{v'+1} = \cfrac{1}{1+\cfrac{1}{v'}} 
= \cfrac{1}{1+\cfrac{1}{[q'_0, q'_1,\ldots, q'_k +  z] }} \\
& = \cfrac{1}{1 + [q'_1,\ldots, q'_k +  z] } 
= [0, 1 + q'_1,\ldots, q'_k +  z] \\
& =  [q_0, q_1,q_2, \ldots, q_k +  z] 
\end{align*}
with $q_0 = 0$, $q_1 = 1 + q'_1$, and $q_i = q'_i$ for $i=2,3,\ldots, k$.  
Moreover, $\sum_{i=0}^k q_i = 1 + \sum_{i=0}^k q'_i = n+1$.

The argument is the same when $v' = [q'_0, q'_1,\ldots, q'_k, z] $ is of form~\eqref{Forest:StandardFormB}.
This completes the proof.  
\end{proof}

\section{The forest of linear fractional transformations}
The set of vertices in the rooted infinite binary tree $\mct(z)$ is the set of all 
special positive linear fractional transformations.  
We constructed this tree by applying the generation rule~\eqref{Forest:parent}:
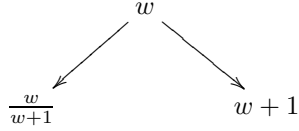
\[
\xymatrix{
& w  \ar[dl]  \ar[dr] & \\
\frac{w }{w +1} &  & w +1 \\
}
\]
to the root $f(z) = z$, and to every successive generation of vertices.
It is a simple observation that this generation rule allows \emph{every} 
positive linear fractional transformation to be the root of an infinite binary tree.
Moreover, if $w = f(z)$, then $w/(w+1) = \ell_1 \circ f(z)$ and $w+1 = r_1\circ f(z)$.
It follows that $\det(w) = \det(w/(w+1) ) = \det(w+1)$.
Thus, if the positive linear fractional transformation $w^*$ is the root 
of the infinite binary tree $\mct(w^*)$, and if $\det(w^*) = D$, 
then every vertex in the tree $\mct(w^*)$ has determinant $D$.  
Because every positive linear fractional transformation is an orphan 
or the descendent of an orphan, it follows that the set of all positive 
linear fractional transformations is a forest of pairwise disjoint rooted infinite binary trees
whose roots are the orphan positive linear fractional transformations.

\bt
Let \mcf\ be the set of positive linear fractional transformations, 
and let $\mco$ be the set of positive linear fractional transformations 
that are orphans.
The infinite forest 
\[
\{ \mct(w^*):w^* \in \mco\}
\]
is a partition of \mcf\ into pairwise disjoint rooted infinite binary trees.

For every nonzero integer $D$, let $\mcf(D)$ be the set of positive linear 
fractional transformations of determinant $D$, 
and let $\mco(D)$ be the set of orphan positive linear fractional transformations 
of determinant $D$.
The set $\mco(D)$ is finite, and the finite forest 
\[
\{ \mct(w^*):w^* \in \mco(D) \}
\]
is a partition of $\mcf(D)$ into pairwise disjoint rooted infinite binary trees.
\et

\begin{proof}
The only statement left to prove is the finiteness of $\mco(D)$.
Let $f(z) = (az+b)/(cz+d)$ be an orphan of determinant $D = ad-bc$.
If $D > 0$, then $a>c$ and $b<d$.  
Equivalently, $a \geq c+1$ and $d \geq b+1$, and so
\[
D = ad-bc \geq (b+1)(c+1)-bc = b+c+1.
\]
There are only finitely many pairs $(b,c)$ of nonnegative integers 
such that $b+c \leq D-1$, and for each such pair there are only finite many 
pairs of nonnegative integers $(a,d)$ such that $D=ad-bc$.
This proves that the number of orphan positive linear fractional 
transformations of determinant $D >0$ is finite.  
Similarly, the number of orphan positive linear fractional 
transformations of determinant $D < 0$ is finite. 
This completes the proof.
\end{proof}

Let $h(D)$ denote the number of orphan positive linear fractional transformations 
of determinant $D$.  
The function
\[
\frac{az+b}{cz+d} \mapsto \frac{cz+d}{az+b}
\]
is a bijection from $\mcf(D)$ to $\mcf(-D)$,
and its restriction to $\mco(D)$ is a bijection from $\mco(D)$ to $\mco(-D)$.
It follows that 
\[
h(D) = h(-D)
\]
for all nonzero integers $D$.
For example, the unique orphan of determinant 1 is $f(z) = z$, 
and the unique orphan of determinant -1 is $f(z) = 1/z$.
Thus, $h(1) = h(-1) = 1$.

There are four orphans of determinant 2:
\[
2z, \qquad \frac{z}{2}, \qquad  \frac{2z}{z+1}, \qquad \frac{z+1}{2}.
\]
and so $h(2) = 2$.  Their reciprocals are the orphans of determinant -2.

The monoid isomorphism from \mcf\ to $GL_2(\N_0)$ defined by
\[
\frac{az+b}{cz+d} \mapsto \bmat a & b \\ c & d \emat
\]
allows us to interchange the languages of positive linear fractional transformations 
and nonnegative matrices, with no gain or loss of generality.
Thus, an \emph{orphan matrix} of determinant $D>0$ is a matrix 
$\bmat a & b \\ c & d \emat$ with $ad-bc=D$ such that  $a>c$ and $b < d$, 
and an \emph{orphan matrix} of determinant $D < 0$ is a matrix 
$\bmat a & b \\ c & d \emat$ with $ad-bc=D$ such that  $a < c$ and $b > d$.

The unique orphan matrix of determinant 1 is 
$I = \bmat 1 & 0 \\ 0 & 1 \emat$, and the unique orphan matrix of determinant -1 is 
$J = \bmat 0 & 1 \\ 1 & 0  \emat$.
The four orphan matrices of determinant 2 are
\[
\bmat
2 & 0 \\
0 & 1
\emat
\qquad
\bmat
1 & 0 \\
0 & 2 
\emat
\qquad
\bmat
2 & 0 \\
1 & 1 
\emat
\qquad
\bmat
1 & 1 \\
0 & 2
\emat.
\]

The seven orphan matrices of determinant 3 are
\[
\bmat  1 & 0 \\ 0 & 3  \emat \quad 
\bmat  3 & 0 \\ 0 & 1  \emat \quad 
\bmat  3 & 0 \\ 1 & 1   \emat \quad 
\bmat   3 & 0 \\ 2 & 1 \emat \quad 
\bmat  1 & 1 \\ 0 & 3  \emat \quad 
\bmat  2 & 1 \\ 1 & 2  \emat \quad 
\bmat  1 & 2 \\ 0 & 3  \emat 
\]
The 13 orphan matrices of determinant 4 are
\[
\bmat  1 & 0 \\ 0 & 4  \emat \quad 
\bmat  2 & 0 \\ 0 & 2 \emat \quad 
\bmat  4 & 0 \\ 0 & 1  \emat \quad 
\bmat   2 & 0 \\ 1 & 2 \emat \quad 
\bmat  4 & 0 \\ 1 & 1   \emat \quad 
\bmat   4 & 0 \\ 2 & 1 \emat \quad 
\bmat   4 & 0 \\ 3 & 1 \emat
\]
\[
\bmat  1 & 1 \\ 0 & 4  \emat  \quad 
\bmat  2 & 1 \\ 0 & 2 \emat \quad 
\bmat  3 & 1 \\ 2 & 2  \emat \quad 
\bmat  1 & 2 \\ 0 & 4 \emat \quad 
\bmat  2 & 2 \\ 1 & 3  \emat \quad 
\bmat 1 & 3 \\ 0 & 4 \emat.
\]

Here is a table of $h(D)$ for $D=1, 2,\ldots, 15$:
\[
\begin{array}{l|ccccccccccccccc}
D      & 1 & 2 & 3 &  4 &   5 &  6 &   7 & 8    & 9 & 10 & 11 & 12  & 13 & 14 & 15 \\ \hline
h(D) & 1 & 4 & 7 & 13 & 15& 26 & 25 & 39& 40 & 54 & 49 & 79 & 63 & 88 & 88 
\end{array}.
\]

\section{A tree grows in a field}
Let $K$ be a field of characteristic 0, and let $z \in K$.  
We would like to construct a rooted infinite  binary tree with root $z$ 
and with vertices in $K$ 
such that every vertex $w$ is the parent of two children: 
the \emph{left child} $w /(w +1)$ 
and the \emph{right child} $w +1$.  
The only obstruction to this construction   
occurs when $w =-1$ for some vertex $w$ in the tree, 
because, in this case, the left child $w/(w+1)$ is undefined.   

Let $K$ be a field of characteristic 0, with multiplicative identity 1.
There is a unique isomorphism from \Q\ to the prime subfield of $K$, 
that is, the subfield of $K$ generated by 1. 
A \emph{negative rational number in $K$} is the image of a negative rational number 
in \Q\ under this isomorphism.  
We shall also denote the prime subfield of $K$ by \Q.  
An element $z \in K$ is called \emph{irrational} if $z\notin \Q$.

\bt     \label{Forest:theorem:creation}
Let $K$ be a field of characteristic 0.
An element in $K$ is the root of a rooted infinite binary tree 
if and only if it is not a negative rational number.
\et

\begin{proof}
If $w\in K\setminus \Q$, then $w+1 \notin \Q$ and $w/(w+1) \notin \Q$.  
Therefore, if $z$ is irrational, then every vertex in the tree 
with root $z$ is irrational.  In particular, $-1$ is not a vertex in this tree.

If $z$ is a positive rational number, then all of its descendants 
are positive rational numbers.  In particular, $-1$ is not a descendant.
If $z=0$, then the left child of 0 is 0 and the right child is 1.  
It follows that the tree of descendants of 0 consists of an infinite 
sequence of 0s, each of which gives birth to a right child 1, 
which is the root of a Calkin-Wilf tree for the positive rational numbers: 
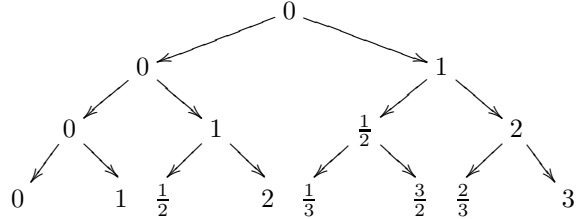
\[
\xymatrix@=0.3cm{
& & & & 0  \ar[dll]  \ar[drr] & & & & \\
& & 0  \ar[dl]  \ar[dr]  &  & & & 1 \ar[dl]  \ar[dr]  & & \\
& 0    \ar[dl]  \ar[dr]  & & 1  \ar[dl]  \ar[dr] & 
 & \frac{1}{2}  \ar[dl]  \ar[dr] &  & 2 \ar[dl]  \ar[dr] & \\
0 &  & 1 \quad \frac{1}{2} &  & 2 \quad  \frac{1}{3} &  & \frac{3}{2} \quad  \frac{2}{3} & & 3  \\
}
\]
Thus, 0 is the root of a rooted infinite binary tree.
It follows that if $z \in K$ is not the root of a rooted infinite binary 
tree, then $z$ is a negative rational number.

Let $z$ be a negative rational number.  
Recall that the \emph{height} \index{height} of a nonzero rational 
number $a/b$ with $a$ and $b$ relatively prime integers 
is $\height(a/b) = \max(|a|, |b|)$.
If $z < 0$ and $\height(z)=1$, then $z= -1$, and so $-1$ is a vertex 
in the tree with root $z$.  
   
Let $h \geq 2$, and suppose that $-1$ is a descendant of every  
negative rational number of height less than $h$.
Let $z$ be a negative rational number of height $h$.  
We can write $z= -a/b$, where $a$ and $b$ are relatively  prime positive integers.  
If $z < -1$, then $1 \leq b < a = \height(z) = h$.  
We have 
\[
 z+1 = -\frac{a-b}{b} < 0
\]
and 
\[
\height(z+1) = \max(a-b,b) < a = h.
\]
The induction hypothesis implies that -1 is a descendent of $z+1$, and so $-1$ is a vertex in the tree with root $z$. 

If $-1 < z < 0$, then $1 \leq  a < b = \height(z) = h$.  We have
\[
\frac{z}{z+1} = -\frac{a}{b-a} < 0
\]
and
\[
\height\left( \frac{z}{z+1}\right)  = \height\left(-\frac{a}{b-a} \right) 
= \max(a,b-a) < b = h.
\]
Again, the induction hypothesis implies that -1 is a descendent of $z/(z+1)$, 
and so $-1$ is a vertex in the tree with root $z$. 
It follows that if $z$ is a negative rational number, then $-1$ is a descendant of $z$, 
and so $z$ is not the root of a rooted infinite binary tree.
This completes the proof.  
\end{proof}

\section{Open problems}
\benum
\item
Find a geometric or algebraic interpretation of the ``class number'' $h(D)$.
Is there a formula to compute $h(D)$?

\item
The subgroup  of $SL_2(\Z)$ generated by 
\[
L_2 =  \bmat 1 & 0 \\ 2 & 1 \emat
\qand
R_2 = \bmat 1 & 2 \\ 0 & 1 \emat
\]
is free of rank 2, and the set $\{L_2, R_2 \}$ freely generates this subgroup
(Sanov~\cite{sano47a}).
Consider the rooted  infinite binary tree that satisfies the generation rule
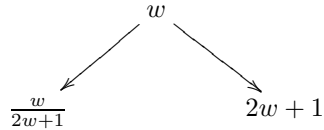
\[
\xymatrix{
& w   \ar[dl] \ar[dr] & \\
\frac{w}{2w+1} & & 2w+1
}
\]
If $w$ is the positive rational number $a/b$, then this generation rule is
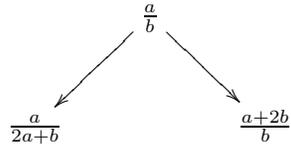
\[
\xymatrix{
& \frac{a}{b}   \ar[dl] \ar[dr] & \\
\frac{a}{2a+b} & & \frac{a+2b}{b}
}
\]
A positive rational number $a/b$ with $\gcd(a,b) = 1$  
is a right child if $a > 2b$ and a left child if $b > 2a$.
The set of positive rational orphans is 
\[
\left\{ \frac{a}{b} : \frac{1}{2} \leq \frac{a}{b}  \leq 2 \right\}.
\]
Every positive rational number generates a rooted infinite binary tree.  
The set of orphan positive rational numbers are the roots of a forest 
of rooted infinite binary binary trees that partition the positive integers.  
Describe the properties of this forest.

\item
Let $u$ and $v$ be  integers such that $u \geq 2$ and $v \geq 2$.  
A standard application of the ping-pong lemma shows 
that the subgroup  of $SL_2(\Z)$ generated by 
\[
L_u =  \bmat 1 & 0 \\ u & 1 \emat
\qand
R_v= \bmat 1 & v \\ 0 & 1 \emat
\]
is free of rank 2, and the set $\{L_u, R_v \}$ freely generates this subgroup
(Lyndon and Schupp~\cite[pp. 167--168]{lynd-schu77}).
Consider the rooted  infinite binary tree that satisfies the generation rule
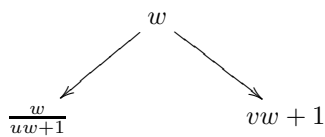
\[
\xymatrix{
& w   \ar[dl] \ar[dr] & \\
\frac{w}{uw+1} & & vw+1
}
\]
If $w$ is the positive rational number $a/b$, then the generation rule is
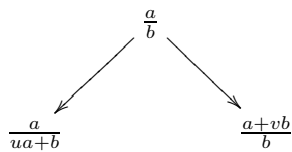
\[
\xymatrix{
& \frac{a}{b}   \ar[dl] \ar[dr] & \\
\frac{a}{ua+b} & & \frac{a+vb}{b}
}
\]
The set of positive rational orphans with respect to the set $\{L_u,R_v\}$  is 
\[
\left\{ \frac{a}{b} : \frac{1}{u} \leq \frac{a}{b}  \leq v \right\}.
\]
Describe the properties of the trees descended from these roots.  

\item
It would be worthwhile to consider trees (not necessarily binary) 
of positive rational numbers and 
of nonnegative matrices constructed from the generators of other free submonoids 
and subgroups of $SL_2(\Z)$ (cf.  Goldberg and Newman~\cite{gold-newm57}
and Bachmuth and Mochizuki~\cite{bach-mozh76}).

\item
Here are some questions about the positive rational number Calkin-Wilf tree.
\benum

\item
Let $a^*/b^*$ and $a/b$ be  positive rational numbers.  
Is there an algorithm to determine if $a/b$ is in the subtree 
of the Calkin-Wilf tree with root $a^*/b^*$?

\item
Is there a continued fraction algorithm to determine the relative depth 
of a vertex $a/b$ is in the subtree 
of the Calkin-Wilf tree with root $a^*/b^*$?

\item
Describe the structure of Calkin-Wilf trees whose roots are negative  rational numbers.  
Does -1 appear infinitely often in a tree with a negative rational root?  

\eenum

\item
Forests of Gaussian numbers are also of interest.

\item
One might also investigate the directed graphs constructed 
with the generation rule~\eqref{Forest:parent} in 
a finite field or in an infinite field of characteristic $p$.

\eenum

\def\cprime{$'$} \def\cprime{$'$} \def\cprime{$'$}
\providecommand{\bysame}{\leavevmode\hbox to3em{\hrulefill}\thinspace}
\providecommand{\MR}{\relax\ifhmode\unskip\space\fi MR }
\providecommand{\MRhref}[2]{%
  \href{http://www.ams.org/mathscinet-getitem?mr=#1}{#2}
}
\providecommand{\href}[2]{#2}

\end{document}